\documentclass{amsart}

\usepackage{amssymb}

\usepackage{url}
\usepackage{color}

\newcommand{\Z}{\mathbb{Z}}

\newcommand{\N}{\Z_{>0}}
\newcommand{\bigoh}{{\mathcal{O}}}

\newtheorem{theorem}{Theorem}[section]

\theoremstyle{definition}
\newtheorem{definition}[theorem]{Definition}

\theoremstyle{remark}

\numberwithin{equation}{section}

\begin{document}

\title{Numerical Verification of the Ternary Goldbach Conjecture up to $8.875\cdot 10^{30}$}


\author{H. A. Helfgott}
\address{Ecole Normale Sup\'erieure, D\'epartement de Math\'ematiques, 45 rue d’Ulm, F-75230 Paris, France}
\email{harald.helfgott@ens.fr}
\author{David J. Platt}
\address{Heilbronn Institute for Mathematical Research, University of Bristol}
\email{dave.platt@bris.ac.uk}
\thanks{We would like to thank the staff of the Direction des Syst\`{e}mes d'Information at Universit\'e Paris VI/VII (Pierre et Marie Curie) and Bill Hart at Warwick University for allowing us to use their computer systems for this 
project. We would also like to thank
 Andrew Booker for suggesting we exploit Proth primes. Travel was funded in part
by the Leverhulme Prize.}



\date{}


\begin{abstract}
We describe a computation that confirms the ternary Goldbach Conjecture up to $8,875,694,145,621,773,516,800,000,000,000$ ($>8.875\cdot 10^{30}$).
\end{abstract}

\maketitle

\section{Introduction}

The ternary Goldbach Conjecture, which has been open since $1742$, is the assertion that every odd number $n>5$ is the sum of three primes. 
The binary, or strong, Goldbach Conjecture states that every even number $>2$ is the sum of two primes; this would imply ternary Goldbach trivially.

Using results from
analytic number theory, specifically Theorems $2$ and $3$ and Table $1$ of
\cite{Ramare2003}, and the truth of the Riemann Hypothesis to $T=3.33\cdot
10^{9}$, Ramar\'e and Saouter showed that ternary Goldbach holds to
at least $8.37\cdot 10^{26}$ (see their Theorem $2.1$). In the first
version of
\cite{Helfgott2012}, the first author used the same argument to show that the verification of the Riemann Hypothesis by the second author \cite{Platt2011} (to $T>3.061\cdot 10^{10}$), Wedeniwski \cite{Wedeniwski2002} (to $T>2.419\cdot 10^{11}$) or Gourdon \cite{Gourdon2010} (to $T>2.4599\cdot 10^{12}$) imply ternary Goldbach holds to $1.23163\cdot 10^{27}$, $6.15697\cdot 10^{28}$ or $5.90698\cdot 10^{29}$ respectively.

In contrast, the most extensive previous explicit verification of ternary Goldbach was by Saouter \cite{Saouter1998} in 1998 who showed that ternary Goldbach holds to $10^{20}$. He in turn relied on a 1993 result of Sinisalo who showed that the binary conjecture holds to $4\cdot 10^{11}$ \cite{Sinisalo1993}. Saouter's approach was to construct a ``ladder'' of primes starting below
$4\cdot 10^{11}$ where each rung of the ladder is no more than $4\cdot 10^{11}$ 
wide. (That is, $p_{i+1}-p_i \leq 4\cdot 10^{11}$.)
Clearly, such a ladder establishes ternary Goldbach up to $4\cdot 10^{11}$ past the last rung and does not rely on any numerical verification of the Riemann Hypothesis. 

Since Saouter's result, the binary Goldbach Conjecture 
has been numerically verified to $4\cdot 10^{18}$ by Oliviera e Silva, 
Siegfried Herzog and Silvio Pardi \cite{Oliviera2013a}. Thus we can immediately gain an advantage of a factor of $10^7$. In addition, there have been significant improvements in hardware capabilities in the intervening $15$ years. We will now describe a computation that exploits both advances.

\section{Proth Primes}

Proving a number of general form is prime can be computationally expensive. 
For example, testing a number of size $N$ using Elliptic Curve Primality Proving (ECPP) has time complexity (heuristically) 
$\bigoh(\log^{4+\epsilon} N)$. Fortunately, there exist primes of special forms that are much easier to test. We use Proth primes\footnote{Saouter exploited numbers of the form $N=2^{22}\cdot R+1$ with $N < 10^{20}$, $N\equiv 3 \bmod 10$ and $R$ odd.}

\begin{definition}
A Proth number is of the form $k\cdot 2^n+1$ with $k,n\in\N$ and $k<2^n$.
\end{definition}

\begin{definition}
A Proth prime is a Proth number that is also prime.
\end{definition}

\begin{theorem}[Proth's Theorem]
A Proth number $N=k\cdot 2^n+1$ is prime if for some integer $a$ with
\begin{equation}\label{eq:schlummer}
\left(\frac{a}{N}\right)=-1
\end{equation}
we have
\begin{equation}\label{eq:siesta}
a^\frac{N-1}{2}\equiv -1 \mod N.
\end{equation}
\begin{proof}
See, for example, $\S$ $5.3.3$ of \cite{Fine2007}.
\end{proof}
\end{theorem}
(Conversely, if $N$ is prime, then (\ref{eq:siesta}) holds for every $a$
satisfying (\ref{eq:schlummer}).)

This suggests the following algorithm for testing a Proth number for primality.
\begin{enumerate}
\item Let $n,k\in\N$ with $k<2^n$. Fix a small prime bound $B$.
\item Set $N\leftarrow k\cdot 2^n+1$.
\item Set $p\leftarrow 2$.
\item While $p<=B$ if $\left(\frac{p}{N}\right)=-1$, then break, else
set $p$ to the first prime after $p$.
\item If $p>B$ (we ran out of small primes) return false.
\item Return true if $p^\frac{N-1}{2} \equiv -1 \mod N$ and 
false if $p^\frac{N-1}{2} \not\equiv -1 \mod N$.
\end{enumerate}

The modular exponentiation at step $6$ requires $\bigoh(\log N)$ multiplications. Computing the Jacobi symbol at step $4$ has (since $a$ is bounded) time complexity $\bigoh(\log(N))$; we do it $\bigoh(1)$ times. 
The algorithm may return a false negative if none of the primes $\leq B$ are quadratic non-residues but $N$ happens to be prime. However, the algorithm never returns a false positive.

\section{The Prime Ladder Algorithm}

We can now formulate the algorithm to construct our prime ladder. 
\begin{enumerate}
\item Fix $n,\Delta\in\N$ with $\Delta>2^n$. Fix a prime $N_0$.
\item Let $k_0$ be the largest integer such that $k_0\cdot 2^n+1<N_0$ and $k_1$ the largest integer such that $k_1\cdot 2^n+1<N_0+\Delta$.
\item While $k_1>k_0$ if $k_1\cdot 2^n+1$ is a Proth prime then break else $k_1\leftarrow k_1-1$.
\item If $k_1=k_0$ (no Proth prime was found) set $N_0$ to the largest 
prime (of general form) less than $N_0+\Delta$; else set $N_0\leftarrow k_1\cdot2^n+1$.
\item Repeat from $2$.
\end{enumerate}

\section{Implementation}

We implemented the above algorithm in the C programming language using the GMP package \cite{GMP2013} to manage integers larger than $64$ bits. To find general primes when no suitable Proth prime could be located, we used Pari's ``precprime'' library routine \cite{Batut2000}. Since ``precprime'' returns a probable prime, we checked each instance with Pari's ``isprime''\footnote{In fact, the numbers returned by ``precprime'' were always prime}. 

We initially hoped to reach about $10^{31}$; this
 required that $n$ be at least $52$. We set $B$ to $29$ so we had $10$ potential candidates for $a$. As each one has about a $50\%$ chance of success, we failed to find a suitable $a$ in about $0.1\%$ of cases.

To improve run time at the expense of space, we test the Proth numbers relating to a given range of $k$ for divisibility by ``small'' primes. We do this by sieving on the arithmetic progression $k\cdot 2^n+1$. We aimed to exploit the multi-core architecture of modern CPUs so we fixed the maximum width of the interval so that there was sufficient memory available for each core to run our code in parallel. We used intervals of width $2^{54}\cdot 10^9$ each containing about $4\cdot 10^9$ candidate Proth numbers. Experiments indicated that sieving for all prime divisors $<16,000$ was an appropriate compromise\footnote{For comparison, Saouter sieved with the first $10$ primes.}.

 It took about $270$ seconds to sieve such an interval and construct a ladder of primes starting within $2\cdot 10^{18}$ of the left end point, ending within $2\cdot 10^{18}$ of the right end point and where no gap between primes exceeded $4\cdot 10^{18}$. This included writing out the values of $k$ and $a$ for every Proth prime found, and storing any general primes found using Pari.

The values of $k$ and $a$ stored were used to check that the relevant Proth number was indeed prime (using Proth's theorem again) and that they were sufficiently closely spaced to establish the required ladder. This check program was written in C++ using the Class Library for Numbers package \cite{cln}, configured not to use GMP so that the primality and spacing of the Proth numbers were confirmed independantly of GMP. As expected, no errors in either package were evident. The check program took a further $40$ seconds for each interval running on a single core and the data file was then deleted. Again the memory demands of the check program were sufficiently modest that all the cores in a CPU could be utilised simultaneously.

In total, we checked $492,700$ ranges of width $2^{54}\cdot 10^9$,
requiring about $40,000$ core hours. We used three resources, a $48$ core,
AMD Magny Cours cluster at Warwick University and two PowerPC clusters at
Universit\'e de Paris VI/VII (UPMC - DSI - P\^{o}le Calcul).
 All three benefit from ECC memory.

The $130,917$ primes of general form we had to find using Pari were 
independently checked using Fran\c{c}ois Morain's ECPP program 
\cite{Morain1992}. This took a trivial amount of computer time and again 
revealed no discrepancies.

\begin{theorem}
Every odd number between $7$ and $T$, where
\begin{equation*}
T=8,875,694,145,621,773,516,800,000,000,000
\end{equation*}
can be written as the sum of three primes.
\end{theorem}

We note that this is factor of $>8\cdot 10^{10}$ higher than the 1998 result. As already observed, an improvement of $\times 10^7$ can be attributed to starting with the stronger statement on binary Goldbach and our computation consumed about $10$ times as much CPU time. The remaining factor of $800$ is at least partially down to the increased clock speed and word length of modern CPUs. Increases in the speed and capacity of memory also allow us to sieve more efficiently and this will also have contributed to the higher throughput we observed.

\section{Going Further}

Clearly, continuing this computation with the current parameters, it is only a matter of CPU cycles to reach $(2^{52}-1)\cdot 2^{52}+1\approx 2\cdot 10^{31}$. However, the first author's proof that ternary Goldbach holds unconditionally above $10^{29}$ \cite{Helfgott2012}\cite{Helfgott2013} renders any such computation nugatory.


\bibliographystyle{amsplain}
\providecommand{\bysame}{\leavevmode\hbox to3em{\hrulefill}\thinspace}
\providecommand{\MR}{\relax\ifhmode\unskip\space\fi MR }
\providecommand{\MRhref}[2]{%
  \href{http://www.ams.org/mathscinet-getitem?mr=#1}{#2}
}
\providecommand{\href}[2]{#2}

\end{document}